\documentclass[12pt, reqno]{article}
\usepackage{amsmath, amssymb, amscd, eucal}
\usepackage{latexsym}
\usepackage{graphicx}
\usepackage{psfrag}

\newtheorem{thm}{Theorem}[section]

\newtheorem{pro}{Proposition}[section]

\newenvironment{prf}{\medskip 
\noindent {\em Proof.}}{\hfill $\square$ \\}

\let\bbar=\overline

\def\ie{i.\ e.\ }

\def\ra{\rightarrow}

\def\({\left(}                   
\def\){\right)}

\def\[{\left[}                   
\def\]{\right]}

\def\operatorname#1{{\rm #1\,}}                                
\def\op{\operatorname}

\def\tms{\times}                 
                  
\def\id{\equiv}

\def\grd{\nabla}         
\def\del{\partial}       
\def\lap{\triangle}

\def\e{\epsilon}         
            
\def\th{\theta}          
           
\def\lam{\lambda}        
\def\m{\mu}

\def\p{\pi}       
\def\r{\rho}

\def\s{\sigma}    
\def\ph{\phi}

\def\G{\Gamma}

\def\Lam{\Lambda} 
       
\def\S{\Sigma}  

\def\Ph{\Phi}     
     
\def\Om{\Omega}

\def\R{\mathbb R}
\def\C{\mathbb C}
\def\H{\mathbb H}

\newcommand{\bX}{\overline X}

\newcommand{\ch}{\mathcal C}
\newcommand{\har}{\hat{\r}}
\title 
{Holography and the Geometry of Certain 
Convex Cocompact Hyperbolic 3-Manifolds}
\date{September, 2002}
\author   {Xiaodong Wang
\thanks{Department of Mathematics, MIT, Cambridge, MA 02139. 
       \tt{Email:xwang@math.mit.edu} }}

\begin{document}
\maketitle
\bibliographystyle{amsalpha}
\pagestyle{plain}
\section{Introduction}
Applying the idea of AdS/CFT correspondence, Krasnov \cite{kr} studied
a class of convex cocompact hyperbolic 3-manifolds. In physics literature
they are known as Euclidean BTZ black holes. Mathematically they can be
described as $\H^3/\G$, where $\G\subset PSL(2,\C)$ is a Schottky group.
His main result, roughly speaking,
identifies the renormalized volume of such a manifold
with the action for the
Liouville theory on the conformal infinity. See Takhtajan and Teo \cite{TT}
for a rigorous proof and related topics.

This is a nice result establishing another holography correspondence.
But the Liouville theory is not yet fully established and the action
which was proposed by Takhtajan and Zograf \cite{ZT} is quite complicated,
so it is desirable to clarify the meaning of the renormalized volume
in a more geometric and transparent way. This question was first 
raised by Manin and Marcolli \cite{MM} and they speculated that the 
renormalized volume could be calculated through the volume of the convex
core of the bulk space based on an explicit example and a recent result
by Brock \cite{Brock} in a different but related situation.

In this paper we try to compute the renormalized volume in terms of geometric data.
As the first step, we compute the renormalized volume using a different
normalization which geometrically is very natural as it uses the distance
function to the convex core. The result is very simple
and geometric.
We first describe the result in the Fuchsian case.
Let $\G\subset PSL(2,\R)$ be a Fuchsian Schottky group with $2g$ generators. 
Let $\Om(\G)\subset S^2$ be its ordinary set.  
In physics, 
$X=\H^3/\G$ is known as the Euclidean version of a non-rotating BTZ black
hole. Mathematically $X$ is a convex cocompact hyperbolic 3-manifold
with the conformal infinity $\S=\Om(\G)/\G$ which is a compact Riemann
surface of genus $g$.
$X$ has a totally geodesic surface $M=\H^2/\G$. In analogy with
general relativity, we can view $M$ as the $t=0$ slice and $X$ as obtained
by evolving $M$. The noncompact hyperbolic surface $M$ has a number of ends.
For each end $E_i$ there is an ``event horizon'' $C_i$, which is a closed
geodesic. Outside $C_i$ the geometry is totally understood as we know
the end is $\R^+\tms S^1$ with the metric $dt^2+\cosh^2(t)d\th^2$, where
$\th$ is periodic with period $L_i$=length($C_i$). The region inside
all the event horizons is precisely the convex core $\ch$ 
and may have some wormholes. This is a compact hyperbolic surface with
a totally geodesic boundary consisting of the closed geodesics $C_i$. 
By Gauss-Bonnet theorem and some topological consideration the area of 
$\ch$ is a topological invariant ($=2\pi(g-1)$) and does not capture 
the geometric information. The result turns out to be given in terms of
$L_i$.

\begin{thm}
The renormalized volume of $X$ by the distance function to the convex
core is given by
\begin{equation}
V=-\frac{\pi}{2}\sum_iL_i.
\end{equation}
\end{thm}

Unfortunately the normalization we use does not give rise to the hyperbolic
metric on the boundary, so what we compute is not the canonical 
renormalized volume $V_c$, which according to Krasnov \cite{kr} is the
Liouville action on the boundary. We will comment on the difference.

In the general case when $\G$ is non-Fuchsian, the picture is much more
complicated. The convex core $\ch$ is then a compact domain whose boundary
is a ``pleated'' hyperbolic surface, \ie a hyperbolic surface with
a measured geodesic lamination. It has only finitely many closed leaves $C_i$
with non-zero bending angle $\th_i$. Let $L_i$ be the length of $C_i$.
The result is

\begin{thm}
The renormalized volume of $X$ by the distance function to the convex core
on $\S$ is given by
\begin{equation}
V=\op{Vol}(\ch)-\frac{1}{2}\sum_i(\pi-\th_i)L_i.
\end{equation}
\end{thm}
  
\bigskip

Again, we emphasize this is not the canonical renormalized volume. We hope
these results are helpful in understanding the geometry and the canonical
renormalized volume. Our results reduce it to a problem on the convex
core. Understanding the difference between the two normalizations raises
a lot of interesting questions.

The paper is organized as follows.
In Section 2, we summarize briefly some background
knowledge. In Section 3, we discuss the geometry of Schottky group.
The computation for the Fuchsian case is then carried out in Section
4. In the last section we discuss the non-Fuchsian case.

\section{Renormalized volume and conformal anomaly}
We first give a very brief introduction to conformally cocompact
Einstein manifolds, the mathematical framework for AdS/CFT correspondence.
Let $\bX$ be a compact manifold of $n+1$ dimensions with boundary $\S$. 
If $r$ is 
a smooth function on $\bX$ with a first order zero on the boundary of 
$\bX$, positive on $X$, then $r $ is called a defining function. 
A Riemannian metric $g$ on $X=\text{Int}\bX$ is called conformally
compact if for any defining function $r $, $\bbar g=r^2g $ extends 
as a smooth metric on $\bX$. 
The restriction of $\bbar g $ to $\S$ gives a metric on $\S$. This metric 
changes
by a conformal factor if the defining function is changed, so $\S$ has a 
well-defined conformal structure $c$. We call $(\S,c)$ 
the conformal infinity of $(X, g)$. If $g$ satisfies the
Einstein equation $\op{Ric}(g)+ng=0$ we say $(X, g)$ is a conformally
compact Einstein manifold.

A metric $h\in c$ on $\S$ determines a unique ``good'' defining function $r$ in
a collar neighborhood of $\S$ such that 
$$g=r^{-2}(dr^2+h_r),$$
where $h_r$ is an $r-$dependent family of metrics on $\S$ with $h_r|_{r=0}=h$.
By the Einstein equation the expansion of $h_r$ is of the following form
(see e. g. Graham\cite{Gra}). For $n$ odd, 
\begin{equation}\label{eodd}
h_{r}=
h_{(0)}+h_{(2)}{r}^2+(even\  powers)+h_{(n-1)}{r}^{n-1}+h_{(n)}{r}^n+\ldots, 
\end{equation}
where the $h_{(j)}$ are tensors on $\S$, and $h_{(n)}$ is trace-free with
respect to $h$. The tensors $h_{(j)}$ for $j\leq n-1$ are locally formally
determined by the metric $h$, but $h_{(n)}$ is formally undetermined.

For $n$ even the analogous expansion is
\begin{equation}\label{even} 
h_{r}=
h_{(0)}+h_{(2)}{r}^2+(even\  powers)+k{r}^n\log r+h_{(n)}{r}^n+\ldots,
\end{equation}
where the $h_{(j)}$ are locally determined for $j\leq n-2$, $k$ is
locally determined and trace-free, but $h_{(n)}$ is formally undetermined.

Consider now the asymptotics of $\op{Vol}(\{r>\e\})$ as $\e\ra 0$.
By the above expansions for $h_r$ we obtain for $n$ odd
\begin{equation}
\op{Vol}(\{r>\e\})=c_0\e^{-n}+c_2\e^{-n+2}+\op{odd\  powers}+c_{n-1}\s^{-1}
+V+o(1)
\end{equation}
and for $n$ even
\begin{equation}
\op{Vol}(\{r>\e\})=c_0\e^{-n}+c_2\e^{-n+2}+\op{even\  powers} 
+c_{n-2}\s^{-2}
-L\log\e +V+o(1).
\end{equation}
The constant term $V$ is called the renormalized volume, which a-priori
depends on the choice $h$ in the conformal class $c$ on $M$.

Actually for $n$ odd it is not difficult to show that $V$ is independent
of the choice of $h$ and thus defines an absolute invariant of the
conformal compact Einstein manifold $(X,g)$. But this is not so if
$n$ is even. 

In dimension $2+1$, conformally compact hyperbolic Einstein manifolds are
precisely the so called convex cocompact hyperbolic 3-manifolds, objects
which have been much studied by geometers since Poincare. 

Let $\G\subset PSL(2,\C)$ be a torsion-free discrete subgroup such that
$X=\H^3/\G$ is noncompact. Let $\Lam(\G)\subset S^2$ be the limit set and
$\Om(\G)$ its complement. The convex core $\ch$ of $M$ is the closed set
$CH(\G)/\G$, where $CH(\G)$ is the convex hull of $\Lam(\G)$ 
in $\H^3$. It is easy to see that $\ch$ is a deformation retract of $X$.
$X$ is convex cocompact if $\ch$ is compact. In this case $X$ is conformally
compact and the 
conformal infinity is the compact Riemannian surface $\S=\Om(\G)/\G$.

Let $h$ be a metric on $\S$ compatible with the conformal structure and
$V_h$ the corresponding renormalized volume.
For another metric 
$\hat{h}=e^{2u}h$, it can be shown
\begin{equation}\label{diff}
V_{\hat{h}}=V_h-\frac{1}{4}\int_S\(|\grd u|^2+Ru\)d\m_h.
\end{equation}
If we know the renormalized volume for one metric $h$ then
the above can be used to calculate the renormalized volume with
respect any other metric $\hat{h}$ in the same conformal class.

If $\S$ has genus $g>1$, then there is a canonical choice of $h$, namely
the hyperbolic metric. 

\begin{pro}
Let $X$ be a convex compact hyperbolic 3-manifold with the conformal
infinity a compact Riemann surface $\S$ of genus $g>1$. Let $h$ be 
the hyperbolic metric on $\S$. Then for any metric $\hat{h}=e^{2u}h$ 
with $Area(\S, \hat{h})=Area(\S, h)=4\pi (g-1)$, we have
\begin{equation}
V_{\hat{h}}\leq V_{h}
\end{equation}
and the identity holds iff $u\id 1$.
\end{pro}

\begin{prf}
We want to show that on the hyperbolic surface $(\S,h)$
\begin{equation*}
E(u)=\int_S\(|\grd u|^2-2u\)dv_h\geq 0
\end{equation*}
for any function $u$ with $\int_{\S}e^{2u}dv_h=1$. By the convexity of the
exponential function we have
$$\int_{\S}u\leq \log \int_{\S}e^{2u}dv_h=0$$
and hence $E(u)\geq 0$. It is obvious that $E(u)=0$ iff $u\id 1$.
\end{prf}

From the above discussion it is natural to expect that the renormalized volume 
calculated by taking
the hyperbolic metric on the conformal infinity can be expressed in
terms of geometric invariants.



\section{The geometry of Schottky 3-manifolds}
We use the ball model $B^3$ for the hyperbolic 3-space and denote its
isometry group by $M(B^3)$, the M\"obius transformations preserving $B^3$. 
A Schottky polyhedron
in $B^3$ is a convex polyhedron $P$ with an even number of sides 
such that no two sides 
of $P$ meet at infinity. Let $\Ph$ be a $M(B^3)$-side-pairing for a Schottky
polyhedron $P$, with $2g$ sides, such that no side of $P$ is paired to 
itself. The group $\G$ generated by $\Ph$ is called a classical Schottky
group of genus $g$. It is a torsion free discrete subgroup of $M(B^3)$
and has $P$ as a fundamental domain. Let $\Om(\G)\subset S^2$ be the 
ordinary set. It is easy to see that $\S=\Om(\G)/G$ is a compact 
Riemann surface of genus $g$. $X=B^3/\G$ is a convex cocompact hyperbolic 
3-manifold with $\S$ as its conformal infinity. Topologically $X$ is
handle body of genus $g$. For details we refer to the book \cite{rat}.

Now we focus on the special case when the Schottky group 
$\G\subset PSL(2,\R)$. Then $X=\H^3/\G$ contains a totally geodesic surface
$M=\H^2/\G$. By considering the exponential map on the normal bundle of 
$M\subset X$, we can write $X=\R\tms M$ such that the metric on $X$ takes
the form
\begin{equation}\label{honc}
g=dr^2+\cosh^2(r)h,
\end{equation}
where $h$ is the hyperbolic metric on $M$. This explicit description
will make the  computation very transparent.

The surface $M$ is noncompact with a finite number of ends. The genus of
the surface and the number of ends depend on the Schottky group $\G$. The 
following figures, which appear in both \cite{kr} and \cite{MM},  show 
two such surfaces. They have the same polyhedron but the side-pairings
are different, and consequently the resulting surfaces are topologically 
different.

\begin{center}
\includegraphics{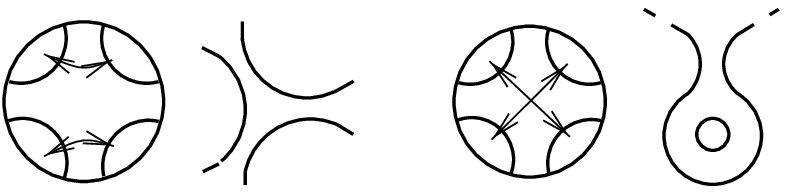}
\end{center}

Suppose $M$ has genus $k$ and $e$ ends. It is easy to get $b_1(M)=k-e+1$.
On the other hand we know $X$ is a handlebody of genus $g$ and hence
$b_1(X)=g$. Since
$M$ is a deformation retract of $X$, their Betti numbers are equal, \ie
\begin{equation}
g=k-e+1.
\end{equation}

For each end there is an outermost closed geodesic $C_i$ and we denote $E_i$
the part outside of $C_i$. Topologically $E_i$ is a cylinder. 
We introduce coordinates $(t,\th)$ on $E_i$ where $t>0$ is the distance
to the boundary $C_i$ and $\th$ is an arc-length parameter on 
$C_i$ and is periodic with period $L_i$=length of $C_i$ . 
Then $E_i=\R^+\tms S^1$ with the metric
\begin{equation}\label{hei}
h=dt^2+\cosh^2(t)d\th^2.
\end{equation}

Cutting off all the ends $E_1,\ldots, E_e$ along these closed 
geodesics $C_1,\ldots,
C_e$, we get a compact hyperbolic surface $\ch$ with totally geodesic
boundary. This is precisely the convex core of both $M$ and $X$.
Therefore we obtain the following decomposition 
\begin{equation}\label{decm}
M=\ch\cup\sqcup_i E_i.
\end{equation}

\section{Computation in the Fuchsian case}
By (\ref{decm}) we also obtain a decomposition for $X$
\begin{equation}
X=(\R\tms\ch)\cup\sqcup_i(\R\tms E_i).
\end{equation}
We define a function $f:X=\R\tms M\ra \R^+$ as follows. For $x\in \ch$
$f(r,x)=|r|$. For $x\in E_i$ we use the coordinates $(t,\th)$ on $E_i$ 
described in
last section and define $f(r,t,\th)>0$ such that $\cosh f=\cosh r\cosh t$.
It is easy to see that $f$ is $C^{1,1}$ and piecewise smooth on 
$X-\ch$  and $|\grd f|\id 1$ by  
(\ref{honc}) and (\ref{hei}). Geometrically  $f$ is just the distance
function to the convex core $\ch$.

Outside the convex core $\ch$, $X$ is foliated by the level sets 
$\S_{\lam}=\{f=\lam\}$ for $\lam\in \(0,\infty\)$. We have the following
decomposition
\begin{equation}\label{dec}
\S_{\lam}=\ch^+(\lam)\cup\ch^-(\lam)\cup\sqcup_i T_i(\lam),
\end{equation}
where 
$$\ch^{\pm}(\lam)=\{\pm\lam\}\tms\ch\subset \R\tms M$$
and
$$T_i(\lam)=\{(r,x)\in \R\tms E_i|\cosh r\cosh t=\cosh\lam\}.$$
With the induced metric both $\ch^+(\lam)$ and $\ch^-(\lam)$ are isometric
to $(\ch, \cosh^2\lam h)$.
We compute the induced metric on $T_i(\lam)$ 
\begin{align*}
dr^2+\cosh^2r(dt^2+\cosh^2td\th^2)&=dr^2
+\frac{\cosh^2\lam\sinh^2r}{\cosh^2\lam-\cosh^2r}dr^2+\cosh^2\lam d\th^2 \\
&=\frac{\cosh^2r\sinh^2\lam}{\cosh^2\lam-\cosh^2r}dr^2+\cosh^2\lam d\th^2 \\
&=\cosh^2\lam\(\frac{\sinh^2 \lam}{\cosh^2\lam}d\ph^2+d\th^2\),
\end{align*}
where $\ph\in \[-\pi/2,\pi/2\]$ is the new variable such that 
$\sin\ph=\frac{\sinh r}{\sinh \lam}$. Therefore in the new coordinates
$(\ph,\th)$ 
\begin{equation}
T_i(\lam)=\[-\pi/2,\pi/2\]\tms S^1
\end{equation} 
with the induced metric taking the following form
\begin{equation}
h_{\lam}=\cosh^2\lam\(\frac{\sinh^2 \lam}{\cosh^2\lam}d\ph^2+d\th^2\).
\end{equation}
This is a standard(flat) cylinder.

Having described all the pieces in the decomposition, we obtain the 
entire surface $\S_{\lam}$ be gluing them together as illustrated by
the following picture.

\begin{center}
\begin{tiny}
\psfrag{t1}[][]{$T_1(\lam)$}
\psfrag{t2}{$T_2(\lam)$}
\psfrag{t3}{$T_3(\lam)$}
\psfrag{ch+}{$\ch^+(\lam)$}
\psfrag{ch-}{$\ch^-(\lam)$}
\includegraphics[totalheight=2.5in]{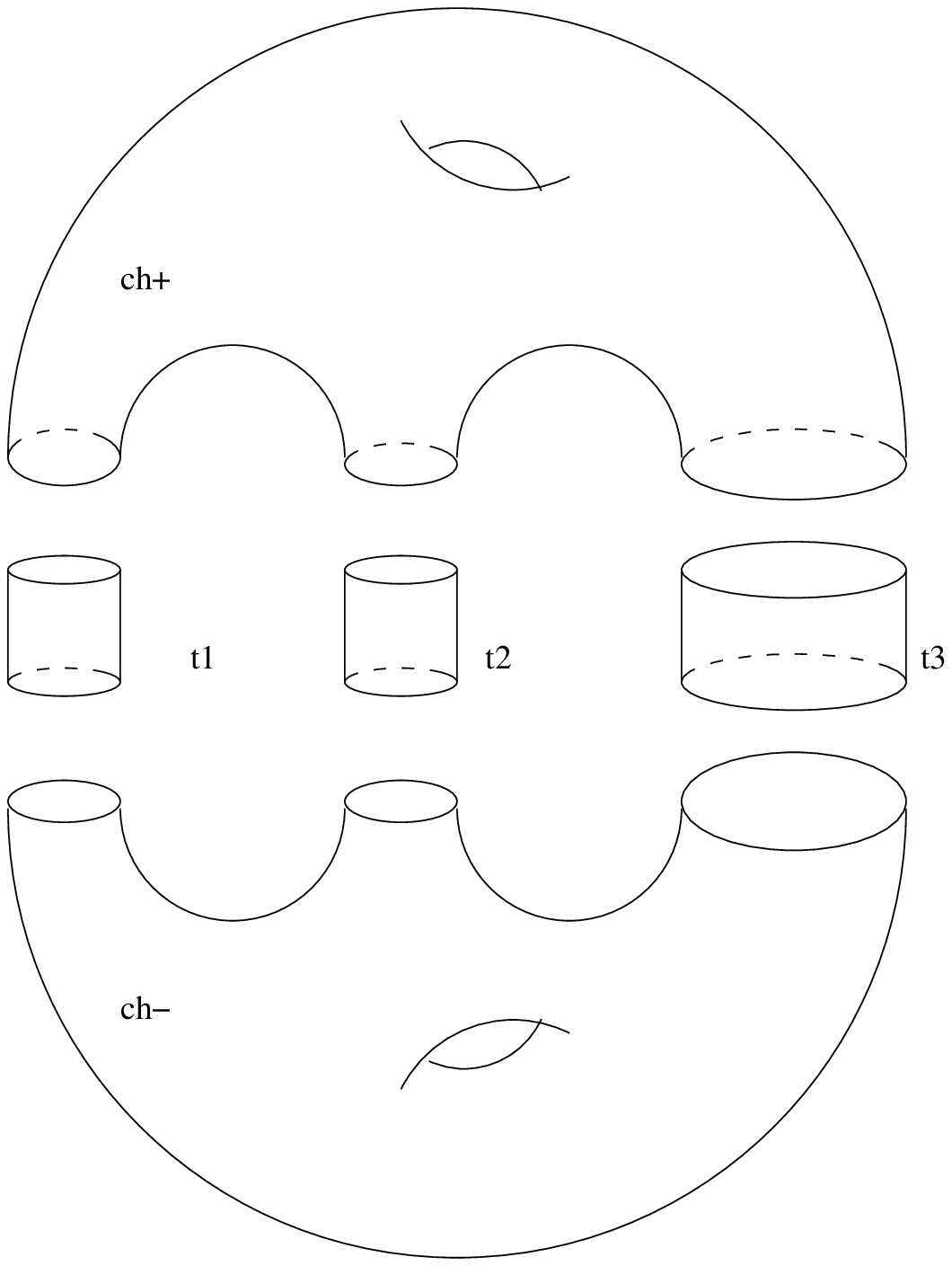}
\end{tiny}
\end{center}

If we scale the metric by dividing the constant factor $\cosh^2\lam$,
the surface $\S_{\lam}$ consists of two copies of the compact hyperbolic 
surface $(\ch, h)$ and a number of cylinders $T_i$ with the base a circle
of length $L_i$ and height $\pi\frac{\sinh \lam}{\cosh \lam}$. As $\lam\ra \infty$ these
height of these cylindrical pieces converges to $\pi$ and we get a closed surface
which consists of two copies of $(\ch, h)$ connected  by these flat cylinders
of height $\pi$ . It is a Riemann surface with a $C^1$ metric, denoted by $h$.
This must be the conformal infinity for $X$.

Let $\har=e^{-f}$. Then we can write  $X-\ch=(0,\infty)\tms \S$ with
the metric $g=\har^{-2}(d\har^2+\hat{h}_{\har}/4)$. We have shown that 
$\hat{h}_{\har}|_{\har=0}$ is the hyperbolic metric on the conformal infinity.

We now compute the renormalized volume. For $\e\in (0,1)$ 
let $X_{\e}=\{(r,x)\in \R\tms M|\har(r,x)> e\}$. We have
\begin{equation}
\op{Vol}(X_{\e})=Vol(\{(r,x)\in X_{\e}|x\in \ch\})+
Vol(\{(r,x)\in X_{\e}|x\in \sqcup_i E_i\}).
\end{equation}
Denote the two summands by $V_1$ and $V_2$. We compute
\begin{equation*}
V_1=\op{Area}(\ch)\int_0^{-\log \e}\cosh^2rdr
=\frac{\pi(g-1)}{4}(\e^{-2}+\frac{\log \e}{2}-\e^2)
\end{equation*}
and
\begin{align*}
V_2&=\sum_i\op{Vol}(\{(r,x)\in \R\tms E_i\cap X_{\e}\}) \\
&=\sum_iL_i\int_{\cosh r\cosh t\leq (\e+\e^{-1})/2}\cosh^2\cosh t drdt \\
&=\frac{\pi}{4}(\e^{-2}-2+\e^2)\sum_iL_i.
\end{align*}
Therefore we obtain 
\begin{equation}
\op{Vol}(X_{\e})=\frac{\pi}{4}(\sum_iL_i+g-1)\e^{-2}+\frac{\pi}{8}(g-1)\log \e
-\frac{\pi}{2}\sum_iL_i+\frac{\pi}{4}(\sum_iL_i-g+1)\e^2.
\end{equation}
The constant term in the above expansion
\begin{equation}
V=-\frac{\pi}{2}\sum_iL_i.
\end{equation}
is then the renormalized volume with respect to $(\S, h)$.

To compute the renormalized volume $V_c$ with respect to the hyperbolic metric $h_0$
on $\S$, we can use formula (\ref{diff}). We write $h_0=e^{2\ph}h$. Then
\begin{align}
\lap \ph+1-e^{2\ph}&=0, \text{on } \ch^+\cup\ch^-, \label{onhy} \\
\lap\ph-e^{2\ph}&=0, \text{ on the flat cylindrical pieces.} \label{onfl}
\end{align}
Then by (\ref{diff}) we have
\begin{equation}
V_c=V-\frac{1}{4}\int_\S\(|\grd \ph|^2+R\ph\)d\m_h
\end{equation}
Note both terms on the right hand side are given on the convex core.
But the second term is very inexplicit as we do not know much about the
$\ph$ which solves (\ref{onhy}) and (\ref{onfl}). It seems difficult to express it
in terms of geometric quantities.
It raises the following general question: Let $S$ be a hyperbolic surface. We get a new 
Riemann surface by cutting it a along a closed
geodesic and then attaching a cylinder of height $t$. How to describe
the hyperbolic metric on the new surface?



\section{The non-Fuchsian case}
We now turn our attention to the general case where the geometry is
much more complicated. The same method works, but the result is less
explicit than the Fuchsian case.
A good reference
for the following discussion is \cite{Epstein}. The original source
is \cite{Thurston}.

For a non-Fuchsian Schottky group $\G\subset PSL(2,\C)$, the limit
set $\Lam(\G)\subset S^2$ is not contained in any great circle. The
convex core $\ch$ is a compact domain in $X$. Its boundary $S=\del\ch$
is a ``pleated'' surface according to Thurston.
With the intrinsic distance $S$ is actually
a hyperbolic surface and how it sits in $X$ is described by a measured
geodesic lamination $K\subset S$. All the geometric information is encoded in
$S$ with this measured geodesic lamination. An important fact is that
$K$ is of measure zero and has only finitely many
closed leaves $C_i$ with nonzero bending angles $\th_i$.

Let $\pi:X-\ch\ra S$ be the nearest point projection. 
Denote $X_{\e}=\{x\in X|d(x, \ch)\leq -\log \e\}$. Then we have the following 
decomposition
\begin{equation}
X_{\e}=\ch\cup \overline{F_{\e}}\cup T_i,
\end{equation}
where,
$$F_{\e}=\{x\in X_{\e}-\ch|\p(x)\in S-K\},$$
and 
$$T_i=\{x\in X_{\e}-\ch|\p(x)\in C_i\}.$$

Let $\har(x)=\exp(-d(x,\ch))$ whcih is $C^{1,1}$ by \cite{Epstein}. 
We claim that this is the defining function
on $X$ that induces $h/4$ on the conformal infinity. To see this we look
at the level set $\S_{\lam}=\{x\in X|d(x,\ch)=\lam\}$. This is a $C^{1,1}$
manifold and we give it the induced metric divided by $\cosh^2(\lam)$.
Since $S-K$ is smooth and totally geodesic in $X$, the piece 
$\{x\in \S_{\lam}|\pi(x)\in S-K\}$ is smooth and hyperbolic. On the other
hand $\{x\in \S_{\lam}|\pi(x)\in \cup C_i\}$ consists of flat pieces which
shrink to disappearance as $\lam \ra \infty$. Therefore $\S_{\lam}$ 
converges to a hyperbolic surface as $\lam\ra \infty$. The detail is 
parallel to the Fuchsian case.

Having shown that  $\har (x)$ is the right defining function, 
the renormalized volume is just the constant term in the asymptotic expansion
as $\e\ra 0$ of
\begin{equation}
\op{Vol}(X_{\e})=\op{Vol}(\ch)+\op{Vol}(F_{\e})+\sum_i\op{Vol}(T_i).
\end{equation}
Since $S-K$ is totally geodesic in $X$, the second piece is very simple. 
The volume of ${F_{\e}}$ is given by
\begin{equation}
V_2=\op{Area}(S)\int_0^{-\log \e}\cosh^2rdr
=\frac{\pi(g-1)}{4}\(\e^{-2}+\frac{\log \e}{2}-\e^2\),
\end{equation}
and there is no contribution to the renormalized volume.
To visualize the pieces $T_i$, we work on the universal covering $\H^3$.
Assume the geodesic $C$ is the $z-$axis in the upper half-space model and the
two bending planes are $y=0$ and $y=\tan (\pi-\th)x$. Then we consider the
following region which is the set of points within distance $r$ to the wedge 
and whose nearest point projection to the wedge lies on the $z$ axis
\begin{equation}
\{(x,y,z)\in \H^3|x\geq 0, y\leq \tan(\th-\pi/2)x, 
\sqrt{x^2+y^2+z^2}/z\leq (\e+\e^{-1})/2 \}.
\end{equation}
It is easy to compute the volume
\begin{equation}
V=(\pi-\th)(\e-\e^{-1})^2/4\int_1^{e^L}dz/z=\frac{(\pi-\th)L}{4}(\e+\e^{-2})
-\frac{(\pi-\th)L}{2} .
\end{equation}
Therefore 
\begin{equation}
\op{Vol}(X_{\e})=\op{Vol}(\ch)-\frac{1}{2}\sum_i (\pi-\th_i)L_i
+\frac{1}{4}\sum_i(\pi-\th_i)L_i(\e+\e^{-2}) \\
+\frac{\pi(g-1)}{4}\(\e^{-2}+\frac{\log \e}{2}-\e^2\),
\end{equation}
and this gives 
the renormalized volume as
$$\op{Vol}(\ch)-\frac{1}{2}\sum_i (\pi-\th_i)L_i.$$

It will be intriguing to see what geometric information can be captured 
by the simple procedure of renormalization for other convex cocompact
hyperbolic manifolds. In a sequel to this paper, we will study 
quasi-Fuchsian deformations.

\providecommand{\bysame}{\leavevmode\hbox to3em{\hrulefill}\thinspace}
\providecommand{\MR}{\relax\ifhmode\unskip\space\fi MR }
\providecommand{\MRhref}[2]{%
  \href{http://www.ams.org/mathscinet-getitem?mr=#1}{#2}
}
\providecommand{\href}[2]{#2}


\begin{thebibliography}{MM01}

\bibitem[Bro]{Brock}
J.~Brock, \emph{The weil-peterssen metric and volumes of 3-dimensional
  hyperbolic convex core}, e-Print math.GT/0109050.

\bibitem[Eps87]{Epstein}
D.~B.~A. Epstein (ed.), \emph{Analytical and geometric aspects of hyperbolic
  space}, Cambridge, Cambridge University Press, 1987. \MR{88c:57003}

\bibitem[Gra00]{Gra}
C.~Robin Graham, \emph{Volume and area renormalizations for conformally compact
  {E}instein metrics}, The Proceedings of the 19th Winter School ``Geometry and
  Physics'' (Srn\'\i, 1999), no.~63, 2000, pp.~31--42. \MR{2002c:53073}

\bibitem[Kra00]{kr}
Kirill Krasnov, \emph{Holography and {R}iemann surfaces}, Adv. Theor. Math.
  Phys. \textbf{4} (2000), no.~4, 929--979. \MR{1 867 510}

\bibitem[MM01]{MM}
Yuri~I. Manin and Matilde Marcolli, \emph{Holography principle and arithmetic
  of algebraic curves}, Adv. Theor. Math. Phys. \textbf{5} (2001), no.~3,
  617--650. \MR{1 898 372}

\bibitem[Rat94]{rat}
John~G. Ratcliffe, \emph{Foundations of hyperbolic manifolds}, Springer-Verlag,
  New York, 1994. \MR{95j:57011}
\bibitem[TT]{TT}
Leon A. Takhtajan amd Lee-Peng Teo, \emph{Liouville action and Weil-Petersson  metric
on deformation spaces, global Kleinian reciprocity and holography}, 
e-print math. CV/0204318. 
\bibitem[Thu]{Thurston}
W.~P. Thurston, \emph{The geometry and topology of three-manifolds}, MSRI 1997,
  electronic version of 1980 Princeton notes, available at {\tt
  http://www.msri.org/gt3m/}.

\bibitem[ZT87]{ZT}
P.~G. Zograf and L.~A. Takhtadzhyan, \emph{On the uniformization of {R}iemann
  surfaces and on the {W}eil-{P}etersson metric on the {T}eichm\"uller and
  {S}chottky spaces}, Mat. Sb. (N.S.) \textbf{132(174)} (1987), no.~3,
  304--321, 444. \MR{88i:32031}

\end{thebibliography}
\end{document}